\documentstyle[11pt]{article}
\begin{document}
\renewcommand{\thefootnote}{\fnsymbol{footnote}}
\pagestyle{empty} \setcounter{page}{1}
\renewcommand{\theequation}{\thesection.\arabic{equation}}
\def\be{\begin{equation}}
\def\eqn#1{\be\label{#1}}

\def\bea{\begin{eqnarray}}
\def\eqnn#1{\bea\label{#1}}
\def\eea{\end{eqnarray}}

\newcommand{\eqna}[1]
{\begin{subequations} \label{#1}
\begin{eqnarray}
\def\eena\end{eqnarray}
\end{subequations}}

\def\id{{\bf 1}}
\def\np{\vfil\eject}
\def\nl{\hfil\break}

\renewcommand{\theequation}{\thesection.\arabic{equation}}
\newfont{\twelvemsb}{msbm10 scaled\magstep1}
\newfont{\eightmsb}{msbm8} \newfont{\sixmsb}{msbm6} \newfam\msbfam
\textfont\msbfam=\twelvemsb \scriptfont\msbfam=\eightmsb
\scriptscriptfont\msbfam=\sixmsb \catcode`\@=11
\def\Bbb{\ifmmode\let\next\Bbb@\else \def\next{\errmessage{Use
      \string\Bbb\space only in math mode}}\fi\next}
\def\Bbb@#1{{\Bbb@@{#1}}} \def\Bbb@@#1{\fam\msbfam#1}
\newfont{\twelvegoth}{eufm10 scaled\magstep1}
\newfont{\tengoth}{eufm10} \newfont{\eightgoth}{eufm8}
\newfont{\sixgoth}{eufm6} \newfam\gothfam
\textfont\gothfam=\twelvegoth \scriptfont\gothfam=\eightgoth
\scriptscriptfont\gothfam=\sixgoth \def\frak{\frak@}
\def\frak@#1{{\fam\gothfam{{#1}}}} \def\frak@@#1{\fam\gothfam#1}
\catcode`@=12

\def\CC{{\Bbb C}}
\def\NN{{\Bbb N}}
\def\QQ{{\Bbb Q}}
\def\RR{{\Bbb R}}
\def\ZZ{{\Bbb Z}}
\def\cA{{\cal A}}          \def\cB{{\cal B}}          \def\cC{{\cal C}}
\def\cD{{\cal D}}          \def\cE{{\cal E}}          \def\cF{{\cal F}}
\def\cG{{\cal G}}          \def\cH{{\cal H}}          \def\cI{{\cal I}}
\def\cJ{{\cal J}}          \def\cK{{\cal K}}          \def\cL{{\cal L}}
\def\cM{{\cal M}}          \def\cN{{\cal N}}          \def\cO{{\cal O}}
\def\cP{{\cal P}}          \def\cQ{{\cal Q}}          \def\cR{{\cal R}}
\def\cS{{\cal S}}          \def\cT{{\cal T}}          \def\cU{{\cal U}}
\def\cV{{\cal V}}          \def\cW{{\cal W}}          \def\cX{{\cal X}}
\def\cY{{\cal Y}}          \def\cZ{{\cal Z}}
\def\qed{\hfill \rule{5pt}{5pt}}
\def\arcsinh{\mathop{\rm arcsinh}\nolimits}
\newtheorem{theorem}{Theorem}
\newtheorem{prop}{Proposition}
\newtheorem{conj}{Conjecture}
\newenvironment{result}{\vspace{.2cm} \em}{\vspace{.2cm}}
\pagestyle{plain}
\begin{center}

{\LARGE {\bf ON SUPER-JORDANIAN ${\cal U}_{\sf h}(sl(N|1))$
ALGEBRA}}
\\[0.2cm]

{\small{\it In memory of our friend Professor Daniel Arnaudon}}

\smallskip

{\small B. ABDESSELAM$^{a}$, A. CHAKRABARTI$^{b}$, R.
CHAKRABARTI$^{c}$, A. YANALLAH$^{d}$ and M.B ZAHAF$^{e}$}

\smallskip

{\small {\it $^{a,d,e}$Laboratoire de Physique Quantique de la
Mati\`ere et Mod\'elisations Math\'ematiques (LPQ3M), Centre
Universitaire de Mascara, 29000-Mascara, Alg\'erie}}

\smallskip

{\small {\it $^{a}$Laboratoire de Physique Th\'eorique d'Oran,
Universit\'e d'Oran Es-S\'enia, 31100-Oran, Alg\'erie}}

\smallskip

{\small {\it $^{b}$Centre de Physique Th\'eorique, Ecole
Polytechnique, 91128-Palaiseau cedex, France.}}

\smallskip

{\small {\it $^{c}$Department of Theoretical Physics, University
of Madras, Guindy Campus, Madras 600025, India}}

\end{center}

\begin{abstract}

\noindent A nonlinear realization of the nonstandard
(super-Jordanian) version of ${\cal U}(sl(N|1))$ is given, for all
$N$.

\end{abstract}


\section{Introduction}
\setcounter{equation}{0}

Jordanian and super-Jordanian quantum algebras have been recently
used for several applications in physical problems. For instance,
super-Jordanian ${\cal U}_{\sf h}\left(osp(1|2)\right)$ algebra
has been understood as the $\kappa$-deformation of the symmetry
algebra of the super-conformal mechanics \cite{BLT}. In another
context, integrable deformed Hamiltonian systems have been
introduced \cite{BH1} via Poisson coalgebra associated with
quantized Jordanian ${\cal U}_{\sf h}\left(sl(2)\right)$ algebra.
We believe that fully developed Hopf coalgebraic structure in a
deformed basis for the ${\cal U}_{\sf h}(sl(N|1))$ presented here
will be useful in building and studying similar deformed fermionic
integrable models. Furthermore, using the corepresentation
structure of the function algebra dually related to the universal
enveloping algebra, a general method of constructing
noncommutative (super)spaces has been recently developed
\cite{ACX} in the context of the quantum supergroup
$OSp_{q}(1|2)$. Application of this method to the case of the dual
quantum supergroup $SL_{\sf h}(N|1)$ will lead to new quantum
superspaces inherently containing a {\it dimensional} deformation
parameter. Influenced by these observations here we introduce the
super-Jordanian ${\cal U}_{\sf h}(sl(N|1))$ algebra in a deformed
basis set.

In a series of papers \cite{ACC1,ACC2,ACC3,ACC4}, we have proposed
a new scheme which permits the construction of the nonstandard
version ${\cal U}_{\sf h}({\sf g})$ of an enveloping
(super)algebra ${\cal U}({\sf g})$ by a suitable contraction, from
the corresponding standard ones ${\cal U}_q({\sf g})$. Our method
hinges on obtaining the ${\cal R}_{\sf h}$-matrix, for all
dimensions, of a (super)Jordanian quantum (super)algebra ${\cal
U}_{\sf h}({\sf g})$ from the ${\cal R}_q$-matrix associated to
the standard quantum (super)algebra ${\cal U}_q({\sf g})$ through
a specific transformation ${\sf G}$ (singular in the $q\rightarrow
1$ limit), as follows:
\begin{equation}
{\cal R}_{\sf h}=\lim_{q\rightarrow 1}\left[{\sf G}^{-1}\otimes
{\sf G}^{-1}\right]{\cal R}_{q} \left[{\sf G}\otimes{\sf
G}\right],
\end{equation}
where, for example, ${\sf G}={\sf E}_q\left(\frac{{\sf h}{\hat
e}_{1N}}{q-1}\right)$ for ${\cal U}_q\left(sl(N)\right)$ (${\hat
e}_{1N}$ is the longest positive root generator of ${\cal
U}_q(sl(N))$) and ${\sf G}={\sf E}_{q^2}\left(\frac{{\sf h} {\hat
e}^2}{q^2-1}\right)$ for ${\cal U}_q(osp(1|2))$ (${\hat e}$ is the
fermionic positive simple root generator of ${\cal
U}_q(osp(1|2))$). The deformed exponential map ${\sf E}_q$ is
defined by
\begin{equation}
{\sf
E}_q\left(\eta\right)=\sum_{n=0}^\infty\frac{\left(\eta\right)^n}{\left[n\right]_q!},
\qquad [n]_q=\frac{q^{n}-q^{-n}}{q-q^{-1}},\qquad
[n]_q!=[n]_q\times [n-1]_q!,\qquad [0]_q!=1.
\end{equation}
For the transformed matrix, the singularities, however, cancel
yielding a well-defined construction. This procedure yields a
nonstandard deformation along with a nonlinear map of the ${\sf
h}$-Borel subalgebra on the corresponding classical Borel
subalgebra, which can be artfully extended to the whole
(super)algebra. The Jordanian quantum algebra ${\cal U}_{\sf
h}\left(sl(N)\right)$ arising from the process cited above
corresponds to the classical matrix $r=h_{1N}\wedge e_{1N}$. {\it
Therefore, the universal ${\cal R}_{\sf h}$-matrix of the full
${\cal U}_{\sf h}(sl(N))$ Hopf algebra, obtained, coincides with
the universal ${\cal R}_{\sf h}$-matrix of the ${\cal U}_{\sf
h}(sl(2))$ Hopf subalgebra \cite{BH} associated with the highest
roots}. In the case of ${\cal U}(osp(1|2))$\footnote{The recent
work shows that there exist three distinct bialgebra structure on
$osp(1|2)$ and all of them are coboundary. We therefore have three
distinct quantization of $osp(1|2)$.}, the super-Jordanian quantum
super-algebra ${\cal U}_{\sf h}\left(osp(1|2)\right)$ occurred
from our treatment is associated to the classical matrix
$r=h\wedge e^2-e\wedge e$. The advantages of our technic are: (1)
With an appropriate choice of basis, the Jordanian quantum Hopf
(super)algebra, obtained by our process, can be endowed with a
relatively simpler coalgebraic structure; (2) Our nonlinear map
permits immediate explicit construction of the finite-dimensional
irreducible representations.

Let us just mention that in general, nonstandard quantum algebras
are obtained by applying Drinfeld twist \cite{Drinfeld} to the
corresponding Lie algebras (see \cite{Lit1,Lit2,Lit3,Lit4,
Lit5,Lit6,Lit7} and refs. there in. The twist deformation of
super-algebras was also discussed in the litterature: \cite{ARS,
BLT} (${\cal U}(osp(1|2))$ case), \cite{BLT1} (${\cal
U}(osp(1|4))$ case) and \cite{T} (general super-algebra case). We
will not consider this way here.

The main object of this paper is to present how our contraction
procedure work for ${\cal U}(sl(N|1))$ superalgebra for obtaining
the nonstandard version ${\cal U}_{\sf h}(sl(N|1))$. For
simplicity, we will limit here ourselves to ${\cal U}(sl(2|1))$
and ${\cal U}(sl(3|1))$. The construction of higher dimensional
super-algebras ${\cal U}_{\sf h}(sl(N|1))$ is presented, briefly,
in the end of this paper. The manuscript is organized as follows:
the super-Jordanian quantum super-algebra ${\cal U}_{\sf
h}(sl(2|1))$ is introduced via a nonlinear map and proved to be a
Hopf Algebra. Higher dimensional super-algebras ${\cal U}_{\sf
h}(sl(N|1))$, $N\geq 3$, are presented in sections 3 and 4. We
conclude in section 5.

\section{${\cal U}_{\sf h}(sl(2|1))$: contraction, nonlinear map and
Hopf structure} \setcounter{equation}{0}

Let us just recall the more important points concerning $sl(2|1)$:
Let $A=\left(a_{ij}\right)$ be the $2\times 2$ matrix given by
$a_{11}=2$, $a_{12}=a_{21}=-1$ and $a_{22}=0$. The Lie Hopf
superalgebra ${\cal U}\left(sl(2|1)\right)$ is generated by the
generators $h_i$, $e_i$ and $f_i$, $i=1,2$, where $h_1$, $h_2$,
$e_1$ and $f_1$ are even
($\deg(h_1)=\deg(h_2)=\deg(e_1)=\deg(f_1)=0$), while $e_2$ and
$f_2$, are odd ($\deg(e_2)=\deg(f_2)=1$), and the commutation
relations
\begin{eqnarray}
&&\left[h_i,\;h_j\right]=0,\qquad \qquad
\left[h_i,\;e_j\right]=a_{ij}e_j,\qquad \qquad
\left[h_i,\;f_j\right]=-a_{ij}f_j,\qquad \qquad
\left[e_i,\;f_j\right]=\delta_{ij}h_i,\nonumber\\
&&\left[e_2,\;e_2\right]=\left[f_2,\;f_2\right]=0, \qquad \qquad
\left[e_1,\left[e_1,e_2\right]\right]=\left[f_1,\left[f_1,f_2\right]\right]=0.
\end{eqnarray}
The two last equations are called the Serre relations. The
commutator $[\;,\;]$ is understood as the $\ZZ_2$-graded one:
$[a,\;b]=ab-\left(-\right)^{\deg(a)\deg(b)}ba$. Defining
\begin{equation}
e_3=e_1e_2-e_2e_1,\qquad f_3=f_2f_1-f_1f_2,\end{equation} we
obtain
\begin{eqnarray}
&&\left[e_1,\;e_3\right]=0, \qquad \qquad
\left[f_3,\;f_1\right]=0,\qquad \qquad  \left[e_2,\;e_3\right]=
0, \qquad \qquad \left[f_2,\;f_3\right]=0,\nonumber \\
&&e_3^2=f_3^2=0,\qquad \qquad \left[e_3,\;f_3\right]=h_1+h_2\equiv
h_3,\qquad \qquad \left[f_1,\;e_3\right]=e_2,\qquad \qquad
\hbox{etc}.
\end{eqnarray}
Let us just mention that there is a $\CC$-algebra automorphism
$\phi$ of ${\cal U}(sl(2|1))$ such that
\begin{equation}
\phi:\left(h_1,h_2,h_3,e_1,e_2,e_3,f_1,f_2,f_3\right)\rightarrow
\left(h_1,-h_3,-h_2,e_1,f_3,-f_2,f_1,-e_3,e_2\right).\end{equation}

The quasitriangular quantum Hopf superalgebra ${\cal
U}_q(sl(2|1))$ ($q$ is an arbitrary complex number), by analogy
with ${\cal U}(sl(2|1))$, is generated by six elements ${\hat
h}_i$, ${\hat e}_i$ and ${\hat f}_i$, $i=1,2$, under the relations
\begin{eqnarray}
&&[{\hat h}_i,\;{\hat h}_j]=0,\qquad \qquad [{\hat h}_i,\;{\hat
e}_j]=a_{ij}{\hat e}_j,\qquad \qquad
[{\hat h}_i,\;{\hat f}_j]=-a_{ij}{\hat f}_j,\nonumber\\
&&[{\hat e}_i,\;{\hat f}_j]=\delta_{ij}\frac{q^{{\hat
h}_i}-q^{-{\hat h}_i}}{q-q^{-1}},\qquad
\qquad {\hat e}_2^2={\hat f}_2^2=0,\nonumber\\
&&{\hat e}_1^2{\hat e}_2-\left(q+q^{-1}\right){\hat e}_1{\hat
e}_2{\hat e}_1+ {\hat e}_2{\hat e}_1^2={\hat f}_1^2{\hat
f}_2-\left(q+q^{-1}\right){\hat f}_1{\hat f}_2{\hat f}_1+ {\hat
f}_2{\hat f}_1^2=0.
\end{eqnarray}
All generators are even except for ${\hat e}_2$ and ${\hat f}_2$
which are odd. The coproducts, counits and antipodes are given by
\begin{eqnarray}
&&\Delta\left({\hat e}_i\right)={\hat e}_i\otimes q^{{\hat
h}_i/2}+q^{-{\hat h}_i/2}\otimes {\hat
e}_i,\qquad\qquad\epsilon({\hat e}_i)=0,\qquad\qquad S({\hat
e}_i)=-q^{{\hat h}_i/2}{\hat e}_iq^{-{\hat h}_i/2},
\nonumber\\
&&\Delta\left({\hat f}_i\right)={\hat f}_i\otimes q^{{\hat
h}_i/2}+q^{-{\hat h}_i/2}\otimes {\hat
f}_i,\qquad\qquad\epsilon({\hat f}_i)=0,\qquad\qquad S({\hat
f}_i)=-q^{{\hat h}_i/2}{\hat f}_iq^{-{\hat h}_i/2},
\nonumber\\
&&\Delta\left({\hat h}_i\right)={\hat h}_i\otimes 1+1\otimes {\hat
h}_i,\qquad\qquad\epsilon({\hat h}_i)=0,\qquad\qquad S({\hat
h}_i)=-{\hat h}_i,
\end{eqnarray}
The universal ${\cal R}$-matrix is given in refs. \cite{KT,Y}.
Note that the definition of the Hopf superalgebra differs from
that of the usual Hopf algebra by the supercommutativity of tensor
product, i.e. $\left(a\otimes b\right)\left(c\otimes
d\right)=(-1)^{\deg(b)\deg(c)}\left(ac \otimes bd\right)$. For
later use, we note that the fundamental representation of (2.5) is
spanned by
\begin{equation}
\begin{array}{lll}
{\hat h}_1=\pmatrix{1&0&0\cr 0&-1&0\cr 0&0&0},\qquad &{\hat
e}_1=\pmatrix{0&1&0\cr 0&0&0\cr
0&0&0},\qquad & {\hat f}_1=\pmatrix{0&0&0\cr 1&0&0\cr 0&0&0},\\
{\hat h}_2=\pmatrix{0&0&0\cr 0&1&0\cr 0&0&1},\qquad & {\hat
e}_2=\pmatrix{0&0&0\cr 0&0&1\cr 0&0&0},\qquad & {\hat
f}_2=\pmatrix{0&0&0\cr 0&0&0\cr 0&1&0}.
\end{array}
\end{equation}

\subsection{Contraction Process}

Following \cite{ACC1}, the ${\cal R}_{\sf h}$ (${\sf h}$ is an
arbitrary complex number) matrix of the super-Jordanian quantum
superalgebra ${\cal U}_{\sf h}(sl(2|1))$, for arbitrary
representations in the two tensor product sectors, can be also
obtained from the ${\cal R}_q$-matrix associated with the
Drinfeld-Jimbo quantum superalgebra ${\cal U}_q(sl(2| 1))$ through
a specific contraction. For simplicity and brevity, let us start
with (fundamental irrep.) $\otimes$ (fundamental irrep.). The
${\cal R}_q$-matrix of ${\cal U}_q(sl(2|1))$ superalgebra in the
(fund.) $\otimes$ (fund.) representation reads
\begin{eqnarray}
\left.R_{\sf h}\right|_{\left(fund.\otimes fund.\right)}=\left(
\begin{array}{ccccccccc}
q&0&0&0&0&0&0&0&0\\
\noalign{\medskip}0&1&0&q-q^{-1}&0&0&0&0&0\\
\noalign{\medskip}0&0&1&0&0&0&q-q^{-1}&0&0\\
\noalign{\medskip}0&0&0&1&0&0&0&0&0\\
\noalign{\medskip}0&0&0&0&q&0&0&0&0\\
\noalign{\medskip}0&0&0&0&0&1&0&q-q^{-1}&0\\
\noalign{\medskip}0&0&0&0&0&0&1&0&0\\
\noalign{\medskip}0&0&0&0&0&0&0&1&0\\
\noalign{\medskip}0&0&0&0&0&0&0&0&-q^{-2}
\end{array}\right).
\end{eqnarray}
The ${\cal R}_{\sf h}$-matrix in the $(fund.\otimes fund.)$
representation is obtained, from (2.9), in the following manner:
\begin{eqnarray}
\left.R_{\sf h}\right|_{\left(fund.\otimes fund.\right)}&=&
\lim_{q\rightarrow 1}\left[{\sf E}^{-1}_{q}\left(\frac{{\sf
h}{\hat e}_1}{q-1} \right)_{fund.}\otimes{\sf
E}^{-1}_{q}\left(\frac{{\sf h}{\hat e}_1}{q-1}\right)_{fund.}
\right]R_q\left[{\sf E}_{q}\left(\frac{{\sf h}{\hat
e}_1}{q-1}\right)_{fund.}\otimes {\sf E}_{q}\left(\frac{{\sf
h}{\hat
e}_1}{q-1}\right)_{fund.}\right]\nonumber\\
&=&\left(\begin{array}{ccccccccc}
1&{\sf h}&0&-{\sf h}&{\sf h}^2&0&0&0&0\\
\noalign{\medskip}0&1&0&0&{\sf h}&0&0&0&0\\
\noalign{\medskip}0&0&1&0&0&0&0&0&0\\
\noalign{\medskip}0&0&0&1&-{\sf h}&0&0&0&0\\
\noalign{\medskip}0&0&0&0&1&0&0&0&0\\
\noalign{\medskip}0&0&0&0&0&1&0&0&0\\
\noalign{\medskip}0&0&0&0&0&0&1&0&0\\
\noalign{\medskip}0&0&0&0&0&0&0&1&0\\
\noalign{\medskip}0&0&0&0&0&0&0&0&-1
\end{array}\right).
\end{eqnarray}
Similarly, using a Maple program\footnote{Our program was
performed for (fund.) $\otimes$ (fund.), (fund.) $\otimes$
(vect.), etc.}, we obtain, for (fundamental irrep.) $\otimes$
(arbitrary irrep.), the following expression:
\begin{equation}
L\equiv \left.R_{\sf h}\right|_{\left(fund.\otimes arb.\right)}
=\pmatrix{T& &-{\sf h}H_1+\frac{\sf h}2\left(T-T^{-1}\right)&&0\cr
&&&&\cr 0&&T^{-1}&&0\cr &&&&\cr 0&&0&&(-1)^F},
\end{equation}
where
\begin{eqnarray}
&&H_1=\frac 12\left(T+T^{-1}\right)h_1=\sqrt{1+{\sf h}^2e_1^2}h_1,
\qquad T^{\pm 1}=\pm {\sf h}e_1+\sqrt{1+{\sf h}^2e_1^2}.
\end{eqnarray}
The above $L$ operator allows immediate construction of the full
Hopf structure of the Borel subalgebra of the ${\cal U}_{\sf
h}(sl(2|1))$ algebra via the FRT formalism\footnote{The algebraic
and coalgebraic properties of the Borel subalgebra are
respectively given by $\left.R_{\sf h}\right|_{\left(fund.\otimes
fund.\right)}L_1L_2=L_2L_1\left.R_{\sf
h}\right|_{\left(fund.\otimes fund.\right)}$,
$\Delta(L)=L\dot{\otimes}  L$, $\varepsilon(L)=1$ and
$S(L)=L^{-1}$.}.

\subsection{Nonlinear map and Hopf structure}

Following refs. \cite{ACC1,ACC3}, let us introduce the generator
\begin{equation}
F_1=f_1-\frac{{\sf h}^2}4e_1\left(h_1^2-1\right).
\end{equation}
We then show that
\begin{eqnarray}
&&TT^{-1}=T^{-1}T=1,\qquad [H_1,T^{\pm 1}]=T^{\pm 2}-1,\nonumber\\
&&[T^{\pm 1},F_1]=\pm \frac{\sf h}2\biggl(H_1T^{\pm 1}+T^{\pm
1}H_1\biggr),\qquad [H_1,F_1]=-{1\over
2}\biggl(TF_1+F_1T+T^{-1}F_1+F_1T^{-1}\biggr),
\end{eqnarray}
with the well known coproducts, counits and antipodes \cite{O}
\begin{eqnarray}
&&\Delta(H_1)=H_1\otimes T+T^{-1}\otimes H_1,\qquad\qquad
\Delta(T^{\pm 1})=T^{\pm 1}\otimes T^{\pm 1},\qquad\qquad
\Delta(F_1)=F_1\otimes T+T^{-1}\otimes
F_1,\nonumber\\
&&S(H_1)=-TH_1T^{-1},\qquad S(T^{\pm 1})=T^{\mp 1},\qquad
S(F_1)=-TF_1T^{-1},\nonumber\\
&&\epsilon(H_1)=\epsilon(F_1)=0,\qquad\epsilon(T^{\pm 1})=1.
\end{eqnarray}
This implies that the Ohn's structure follows from the bosonic
generators $\{h_1,e_1,f_1\}$. The algebraic properties (2.11) and
(2.12) exhibits clearly the embedding of ${\cal U}_{\sf h}(sl(2))$
in ${\cal U}_{\sf h}(sl(2|1))$.

To complete now the ${\cal U}_{\sf h}(sl(2|1))$ superalgebra, we
introduce the following ${\sf h}$-deformed fermionic root
generators:
\begin{eqnarray}
&&H_2=h_2-\frac{{\sf h}^2}2e_1^2h_1,\qquad E_2=e_2-\frac{{\sf
h}^2}4e_1e_3\left(2h_1+1\right),\qquad F_2=f_2,\nonumber\\
&&H_3=h_3+\frac{{\sf h}^2}2e_1^2h_1,\qquad E_3=e_3,\qquad
F_3=f_3+\frac{{\sf h}^2}4e_1f_2\left(2h_1+1\right).
\end{eqnarray}
The generators $E_2$, $E_3$, $F_2$ and $F_3$ are odd, while $H_2$
and $H_3$ are even. The expressions (2.11), (2.12) and (2.15)
define a realization of the super-Jordanian subalgebra ${\cal
U}_{\sf h}(sl(2|1))$ with the classical generators via a nonlinear
map ({\it Other invertible maps relating the super-Jordanian and
the classical generators may also be considered}) and permit
immediate explicit construction of the finite-dimensional
irreducible representations of the ${\cal U}_{\sf h}(sl(2|1))$
superalgebra. In the followings we only quote the final results:
\begin{prop}
The nonstandard (super-Jordanian) enveloping superalgebra ${\cal
U}_{\sf h}(sl(2|1))$ is an associative superalgebra over $\CC$
generated by $\{H_1,\;T,\;T^{-1}
,\;F_1,\;H_2,\;E_2,\;F_2,\;H_3,\;E_3,\;F_3\}$ satisfying, along
with (2.15) and (2.17), the commutation relations
\begin{eqnarray}
&&[H_1,\;H_2]=-\frac 14\left(T-T^{-1}\right)^2H_1,\qquad
[H_1,\;H_3]=\frac 14\left(T-T^{-1}\right)^2H_1,\qquad
[H_2,\;H_3]=0,\nonumber\\
&&[H_1,\;E_2]=-\frac 12\left(T+T^{-1}\right)E_2-\frac{\sf
h}2\left(T-T^{-1}\right)E_3H_1
-\frac{\sf h}4\left(T^2-T^{-2}\right)E_3,\nonumber\\
&&[H_1,\;F_3]=-\frac 12\left(T+T^{-1}\right)F_3+\frac{\sf
h}2\left(T-T^{-1}\right)F_2H_1+
\frac{\sf h}4\left(T^2-T^{-2}\right)F_2,\nonumber\\
&&[H_1,\;F_2]=\frac 12\left(T+T^{-1}\right)F_2,\qquad\qquad
[H_1,\;E_3]=\frac 12\left(T+T^{-1}
\right)E_3,\nonumber\\
&&[H_2,\;T^{\pm 1}]=-\frac 14\left(T^{\pm 3}-T^{\mp
1}\right),\qquad\qquad
[H_3,\;T^{\pm 1}]=\frac 14\left(T^{\pm 3}-T^{\mp 1}\right),\nonumber\\
&&[H_2,\;F_1]=\frac 14\left(T+T^{-1}\right)^2F_1-\frac{\sf
h}4\left(T-T^{-1}\right)H_1^2 -\frac{\sf
h}4\left(T^2-T^{-2}\right)H_1-\frac{\sf
h}{16}\left(T^2-T^{-2}\right)\left(T+
T^{-1}\right),\nonumber\\
&&[H_3,\;F_1]=-\frac 14\left(T+T^{-1}\right)^2F_1+\frac{\sf
h}4\left(T-T^{-1}\right)H_1^2 +\frac{\sf
h}4\left(T^2-T^{-2}\right)H_1+\frac{\sf
h}{16}\left(T^2-T^{-2}\right)\left(T+
T^{-1}\right),\nonumber\\
&& [H_2,\;E_2]=\frac{\sf
h}{16}\left(T+T^{-1}\right)\left(T^2-T^{-2}\right)E_3
+\frac 18\left(T-T^{-1}\right)^2E_2,\nonumber\\
&&[H_3,\;F_3]=\frac{\sf
h}{16}\left(T-T^{-1}\right)\left(T^2-T^{-2}\right)F_2
-\frac 18\left(T-T^{-1}\right)^2F_3,\nonumber\\
&& [H_2,\;F_3]=\frac 18\left(T^2+6+T^{-2}\right)F_3-\frac{\sf
h}{16}\left(T^2-T^{-2}\right)
\left(T+T^{-1}\right)F_2,\nonumber\\
&& [H_3,\;E_2]=-\frac 18\left(T^2+6+T^{-2}\right)E_2-\frac{\sf
h}{16}\left(T^2-T^{-2}\right)
\left(T+T^{-1}\right)E_3,\nonumber\\
&& [H_2,\;F_2]=-\frac 18\left(T-T^{-1}\right)^2F_2,\qquad
[H_3,\;E_3]=\frac 18\left(T-T^{-1}\right)^2E_3,
\nonumber\\
&&[H_3,\;F_2]=\frac 18\left(T^2+6+T^{-2}\right)F_2,\qquad
[H_2,\;E_3]=-\frac 18\left(T^2+6+T^{-2}\right)E_3,\nonumber\\
&&[E_2,\;F_2]=H_2-\frac 1{16}\left(T-T^{-1}\right)^2-\frac{\sf
h}4\left(T-T^{-1}\right)E_3F_2,\nonumber\\
&&[E_3,\;F_3]=H_3+\frac 1{16}\left(T-T^{-1}\right)^2+\frac{\sf
h}4\left(T-T^{-1}\right)F_2E_3,\nonumber\\
&&[T^{\pm 1},\;F_2]=[T^{\pm 1},\;E_3]=0,\qquad
F_2^2=E_3^2=0,\qquad
[F_2,\;F_1]=F_3,\qquad [F_1,\;E_3]=E_2,\nonumber\\
&&E_2^2=\frac{\sf h}4\left(T-T^{-1}\right)E_3E_2,\qquad
F_3^2=-\frac{\sf h}4\left(T-T^{-1}\right)F_2F_3,\qquad
[E_2,\;E_3]=[F_2,\;F_3]=0,\nonumber\\
&&[T^{\pm 1},\;E_2]=\pm \frac{\sf h}2\left(T^{\pm
2}+1\right)E_3,\qquad [T^{\pm 1},\;F_3]=\mp\frac{\sf
h}2\left(T^{\pm 2}+1\right)F_2,\qquad
[F_2,\;E_3]=\frac1{2{\sf h}}\left(T-T^{-1}\right),\nonumber\\
&&[E_2,\;F_1]=\frac{\sf h}4\left(T-T^{-1}\right)E_2+\frac{\sf
h}2\left(T-T^{-1}\right)E_3F_1 -\frac{{\sf
h}^2}4E_3H_1^2-\frac{3{\sf h}^2}8\left(T+T^{-1}\right)E_3H_1-
\frac{{\sf h}^2}2E_3\nonumber\\
&&\phantom{[E_2,\;F_1]=}-\frac{15{\sf
h}^2}{64}\left(T-T^{-1}\right)^2E_3,\nonumber\\
&&[F_3,\;F_1]=\frac{\sf h}4\left(T-T^{-1}\right)F_3-\frac{\sf
h}2\left(T-T^{-1}\right)F_2F_1 +\frac{{\sf
h}^2}4F_2H_1^2+\frac{3{\sf h}^2}8\left(T+T^{-1}\right)F_2H_1+
\frac{{\sf h}^2}2F_2\nonumber\\
&&\phantom{[E_2,\;F_1]=}+\frac{15{\sf
h}^2}{64}\left(T-T^{-1}\right)^2F_2,\nonumber\\
&&[F_3,\;E_2]=F_1-\frac{\sf
h}4\left(T-T^{-1}\right)F_2E_2+\frac{\sf h}4\left(T-T^{-1}\right)
E_3F_3-\frac{\sf h}8\left(T-T^{-1}\right)H_1^2-\frac{\sf
h}8\left(T^2-T^{-2}\right)H_1
\nonumber\\
&&\phantom{[F_3,\;E_2]=}-\frac{\sf
h}{16}H_1\left(T^2-T^{-2}\right)- \frac{7{\sf
h}}{128}\left(T-T^{-1}\right)^3.
\end{eqnarray}
The $\ZZ_{2}$-grading in ${\cal U}_{\sf h}(sl(2|1))$ is uniquely
defined by the requirement that the only odd generators are $E_2$,
$F_2$, $E_3$ and $F_3$. It is obvious that as ${\sf h}\rightarrow
0$, we have $\left(E_2,F_2,H_2,E_3,F_3,
H_3\right)\rightarrow\left(e_2,f_2,h_2,e_3,f_3,h_3\right)$.
\end{prop}
\begin{prop}
Let us note that there exist a $\CC$-algebra automorphism of
${\cal U}_{\sf h}(sl(2|1))$ such that
\begin{equation}
\Phi\left(T^{\pm 1},F_1,H_1,E_2,
F_2,H_2,E_3,F_3,H_3\right)\longrightarrow \left(T^{\pm
1},F_1,H_1,F_3,-E_3,-H_3,-F_2,E_2,-H_2\right).
\end{equation}
(For ${\sf h}=0$, this automorphism reduces to (2.5)).
\end{prop}
\begin{prop}
The nonstandard (super-Jordanian) quantum enveloping superalgebra
${\cal U}_{\sf h}(sl(2|1))$ admits a Hopf structure with
coproducts, antipodes and counits determined by (2.15) and
\begin{eqnarray}
&&\Delta\left(E_2\right)=E_2\otimes T^{1/2}+T^{-1/2}\otimes
E_2+\frac{\sf h}4T^{-1}E_3\otimes
\left(T^{-1/2}H_1+H_1T^{-1/2}\right)-\frac{\sf
h}4\left(T^{1/2}H_1+H_1T^{1/2}\right)\otimes
TE_3,\nonumber\\
&&\Delta\left(F_2\right)=F_2\otimes T^{-1/2}+T^{1/2}\otimes
F_2,\qquad \Delta\left(E_3\right)=E_3\otimes
T^{-1/2}+T^{1/2}\otimes
E_3,\nonumber\\
&&\Delta\left(F_3\right)=F_3\otimes T^{1/2}+T^{-1/2}\otimes
F_3-\frac{\sf h}4T^{-1}F_2\otimes
\left(T^{-1/2}H_1+H_1T^{-1/2}\right)+\frac{\sf
h}4\left(T^{1/2}H_1+H_1T^{1/2}\right)\otimes
TF_2,\nonumber\\
&&\Delta\left(H_2\right)=H_2\otimes 1+1\otimes H_2+\frac
14TH_1\otimes\left(1-T^2\right)+
\frac 14\left(1-T^{-2}\right)\otimes T^{-1}H_1,\nonumber\\
&&\Delta\left(H_3\right)=H_3\otimes 1+1\otimes H_3-\frac
14TH_1\otimes\left(1-T^2\right)-
\frac 14\left(1-T^{-2}\right)\otimes T^{-1}H_1,\nonumber\\
&&S\left(E_2\right)=-E_2-\frac{\sf
h}2\left(T+T^{-1}\right)E_3,\qquad
S\left(F_3\right)=-F_3+\frac{\sf h}2\left(T+T^{-1}
\right)F_2,\nonumber\\
&&S\left(F_2\right)=-F_2,\qquad S\left(E_3\right)=-E_3,\nonumber\\
&&S\left(H_2\right)=-H_2+\frac
12\left(T^{-2}-1\right),\qquad\qquad
S\left(H_3\right)=-H_3-\frac 12\left(T^{-2}-1\right),\nonumber\\
&&\epsilon\left(H_2\right)=\epsilon\left(H_3\right)=\epsilon\left(E_2\right)=\epsilon
\left(F_2\right)=\epsilon\left(E_3\right)=\epsilon\left(F_3\right)=0.
\end{eqnarray}
\end{prop}
All the Hopf superalgebra axioms can be verified by direct
calculations. We remark that our coproducts have simpler forms
compared to those given in the literature. This is one main
advantage of our procedure.
\begin{prop}
The universal ${\cal R}_{\sf h}$-matrix of ${\cal U}_{\sf
h}(sl(2|1))$ has the following form:
\begin{eqnarray}
&& {\cal R}_{\sf h}=\exp\biggl(-{\sf h}X_1\otimes
TH_1\biggr)\exp\biggl({\sf h}TH_1\otimes X_1\biggr),
\end{eqnarray}
where $X_1={\sf h}^{-1}\ln T$. The element (2.19) coincides with
the pure ${\cal U}_{\sf h}(sl(2))$ universal ${\cal R}_{\sf
h}$-matrix \cite{BH}.
\end{prop}

\section{${\cal U}(sl(3|1))$: Nonstandard quantization and Hopf
Structure} \setcounter{equation}{0}

The major interest of our approach is that it can be generalized
for obtaining super-Jordanian quantum superalgebras ${\cal U}_{\sf
h}(sl(N|1))$ of higher dimensions. We start here with ${\cal
U}_{\sf h}(sl(3|1))$. In our notations $e_{ij}$ is an $(N+1)\times
(N+1)$ matrix with only the $(i,j)$ matrix element being equal to
1, all other matrix elements are zero. Let $h_{12}=e_{11}-
e_{22}$, $h_{23}=e_{22}-e_{33}$, $h_{34}=e_{33}+e_{44}$, $e_{12}$,
$e_{23}$, $e_{34}$, $e_{21}$, $e_{32}$ and $e_{43}$ be the
standard Chevalley generators of ${\cal U}\left(sl(3|1)\right)$.
The generators $h_{12}$, $h_{23}$, $e_{12}$, $e_{23}$, $e_{21}$,
$e_{32}$, and $h_{34}$ are even, while $e_{34}$ and $e_{43}$ are
odd. The generators corresponding to the other roots, obtained by
action of the Weyl group, are denoted by $e_{13}=[e_{12},e_{23}]$,
$e_{14}=[e_{13}, e_{34}]$, $e_{24}=[e_{23},e_{34}]$,
$e_{31}=[e_{32},e_{21}]$, $e_{41}=[e_{43},e_{31}]$, $e_{42}
=[e_{43},e_{32}]$, $h_{13}=e_{11}-e_{33}\equiv h_{12}+h_{23}$,
$h_{14}=e_{11}+e_{44}\equiv h_{13}+h_{34}$ and
$h_{24}=e_{22}+e_{44}\equiv h_{23}+h_{34}$\footnote{The elements
$\{h_{12},\;h_{23},
\;e_{12},\;e_{23},\;e_{21},\;e_{32},\;e_{13},\;e_{31},\; h_{13}\}$
build here the subalgebra ${\cal U}(sl(3))$ of ${\cal
U}(sl(3|1))$. }. The commutator $[\;,\;]$ is understood as the
$\ZZ_2$-graded one, i.e.
\begin{equation}
\left[e_{ij},\;e_{kl}\right]=\delta_{jk}e_{il}-(-)^{\deg\left(e_{ij}
\right)\deg\left(e_{kl}\right)}\delta_{li}e_{kj}.
\end{equation}
There exist a $\CC$-algebra automorphism $\phi$ of ${\cal
U}(sl(3|1))$ such that
\begin{equation}
\phi\left(e_{12},e_{21},h_{12},e_{23},e_{32},h_{23},e_{34},e_{43},h_{34},\cdots\right)\longrightarrow
\left(e_{23},e_{32},h_{23},e_{12},e_{12},h_{12},e_{41},e_{14},-h_{14},\cdots\right)
\end{equation}

\subsection{The Bosonic part: ${\cal U}_{\sf h}(sl(3))$ subalgebra}
As in the ${\cal U}_{\sf h}(sl(2|1))$ superalgebra, the
super-Jordanian deformation arises here from the bosonic
generators corresponding to the higher root, i.e. $e_{13}$,
$e_{31}$ and $h_{13}$. These generators are deformed as
follows\footnote{Similarly to Ref. \cite{ACC3}, by applying the
contraction process on the $R_q$-matrix in the $(fund.\otimes
arb.)$, associated to ${\cal U}_q(sl(3|1))$, we obtain:
\begin{eqnarray}
&&\left.R_{\sf h}\right|_{\left(fund.\otimes arb.\right)}
=\pmatrix{T& 2{\sf h}T^{-1/2}e_{23}&-\frac{\sf
h}2\left(T+T^{-1}\right)(h_1+h_2)+\frac{\sf
h}2\left(T-T^{-1}\right)&0\cr &&&\cr 0&I&-2{\sf
h}T^{1/2}e_{12}&0\cr &&&\cr 0&0&T&0\cr 0&0&0&(-1)^F}.\nonumber
\end{eqnarray}}:
\begin{eqnarray}
&&T^{\pm 1}=\pm{\sf h}e_{13}+\sqrt{1+{\sf h}^2e_{13}^2},\qquad
H_{13}=\sqrt{1+{\sf h}^2e_{13}^2}h_{13},\qquad
E_{31}=e_{31}-\frac{{\sf h}^2}4e_{13}\left(h_{13}^2-1\right).
\end{eqnarray}

To complete first the ${\cal U}_{\sf h}(sl(3))\subset {\cal
U}_{\sf h}(sl(3|1))$ subalgebra ({\it the bosonic part of ${\cal
U}_{\sf h}(sl(3|1))$}), let us introduce the following ${\sf
h}$-deformed generators:
\begin{eqnarray}
&&H_{12}=h_{12}+\frac{{\sf h}^2}2e_{13}^2h_{13}, \qquad
E_{12}=e_{12},\qquad E_{21}=e_{21}+\frac{{\sf
h}^2}4e_{23}e_{13}\left(2h_{13}+1\right),\nonumber\\
&&H_{23}=h_{23}+\frac{{\sf h}^2}2e_{13}^2h_{13},\qquad
E_{23}=e_{23},\qquad E_{32}=e_{32}-\frac{{\sf
h}^2}4e_{12}e_{13}\left(2h_{13}+1\right),
\end{eqnarray}
where it is obvious that as ${\sf h}\rightarrow 0$, we have
$\left(H_{12},\;E_{12},\;E_{21},\;H_{23},\;E_{23},\;E_{32},\
;H_{13},\;T,\;T^{-1},\;E_{31}\right)\rightarrow (h_{12},\;e_{12},$
$e_{21},\;h_{23},\;e_{23},\;e_{32},\;h_{13},\;1,\;1, \;e_{31})$.
The expressions (3.3) and (3.4) define a realization of the
Jordanian subalgebra ${\cal U}_{\sf h}(sl(3))$ embedded in ${\cal
U}_{\sf h}(sl(3|1))$ with the classical generators via a nonlinear
map. Another map has been considered in \cite{ACC3}. Our
construction leads to the following results:
\begin{prop}
The generating elements $\{H_{12},\;E_{12},\;E_{21},\;H_{23},
\;E_{23},\;E_{32},\;H_{13},\;T,\;T^{-1},\;E_{31}\}$ of the
Jordanian quantum algebra ${\cal U}_{\sf h}(sl(3))$ obey the
following commutations rules:
\begin{eqnarray}
&&TT^{-1}=T^{-1}T=1,\qquad [H_{13},T^{\pm 1}]=T^{\pm 2}-1, \qquad
[T^{\pm 1},E_{31}]=\pm\frac{\sf h}2\biggl(H_{13}T^{\pm 1}+T^{\pm
1}H_{13}\biggr),\nonumber\\
&&[H_{13},E_{31}]=-{1\over
2}\biggl(\left(T+T^{-1}\right)E_{31}+E_{31}\left(T+T^{-1}
\right)\biggr),\nonumber\\
&&[H_{12},\;H_{23}]=0,\qquad [H_{12},\;H_{13}]=-\frac
14\left(T-T^{-1}\right)^2H_{13},\qquad
[H_{23},\;H_{13}]=-\frac 14\left(T-T^{-1}\right)^2H_{13},\nonumber\\
&& [H_{12},\;E_{12}]=2E_{12}+\frac
18\left(T-T^{-1}\right)^2E_{12},\qquad
[H_{12},\;E_{23}]=-E_{23}+\frac
18\left(T-T^{-1}\right)^2E_{23},\nonumber\\
&& [H_{23},\;E_{12}]=-E_{12}+\frac
18\left(T-T^{-1}\right)^2E_{12},\qquad
[H_{23},\;E_{23}]=2E_{23}+\frac
18\left(T-T^{-1}\right)^2E_{23},\nonumber\\
&& [H_{12},\;E_{21}]=-2E_{21}-\frac
18\left(T-T^{-1}\right)^2E_{21}+\frac{\sf h}{16}
\left(T+T^{-1}\right)\left(T^2-T^{-2}\right)E_{23},\nonumber\\
&& [H_{23},\;E_{32}]=-2E_{32}-\frac
18\left(T-T^{-1}\right)^2E_{32}-\frac{\sf h}{16}
\left(T+T^{-1}\right)\left(T^2-T^{-2}\right)E_{12},\nonumber\\
&& [H_{12},\;E_{32}]=E_{32}-\frac
18\left(T-T^{-1}\right)^2E_{32}-\frac{\sf h}{16}
\left(T+T^{-1}\right)\left(T^2-T^{-2}\right)E_{12},\nonumber\\
&& [H_{23},\;E_{21}]=E_{21}-\frac
18\left(T-T^{-1}\right)^2E_{21}+\frac{\sf h}{16}
\left(T+T^{-1}\right)\left(T^2-T^{-2}\right)E_{23},\nonumber\\
&&[H_{13},\;E_{12}]=\frac
12\left(T+T^{-1}\right)E_{12},\qquad\qquad [H_{13},\;E_{23}]=\frac
12\left(T+T^{-1}\right)E_{23},\nonumber\\
&& [H_{13},\;E_{21}]=-\frac
12\left(T+T^{-1}\right)E_{21}+\frac{\sf h}{2}\left(T-T^{-1}\right)
E_{23}H_{13}+\frac{\sf h}4\left(T^2-T^{-2}\right)E_{23},\nonumber\\
&& [H_{13},\;E_{32}]=-\frac
12\left(T+T^{-1}\right)E_{32}-\frac{\sf h}{2}\left(T-T^{-1}\right)
E_{12}H_{13}-\frac{\sf h}4\left(T^2-T^{-2}\right)E_{12},\nonumber\\
&&[E_{21},\;F_{31}]=\frac{\sf
h}4\left(T-T^{-1}\right)E_{21}-\frac{\sf h}2\left(T-T^{-1}\right)
E_{23}E_{31}+\frac{{\sf h}^2}4E_{23}H_{13}^2+\frac{3{\sf
h}^2}8\left(T+T^{-1}\right)E_{23}H_{13}
\nonumber\\
&&\phantom{[E_2,\;F_1]=}+\frac{{\sf h}^2}2E_{23}+\frac{15{\sf
h}^2}{64}\left(T-T^{-1}\right)^2
E_{23},\nonumber\\
&&[E_{32},\;F_{31}]=\frac{\sf
h}4\left(T-T^{-1}\right)E_{32}+\frac{\sf h}2\left(T-T^{-1}\right)
E_{12}E_{31}-\frac{{\sf h}^2}4E_{12}H_{13}^2-\frac{3{\sf
h}^2}8\left(T+T^{-1}\right)E_{12}H_{13}
\nonumber\\
&&\phantom{[E_2,\;F_1]=}-\frac{{\sf h}^2}2E_{12}-\frac{15{\sf
h}^2}{64}\left(T-T^{-1}\right)^2
E_{12},\nonumber\\
&&[H_{12},\;T^{\pm 1}]=-\frac 14\left(T^{\pm 3}-T^{\mp
1}\right),\qquad [H_{23},\;T^{\pm 1}]=-\frac 14\left(T^{\pm
3}-T^{\mp
1}\right),\nonumber\\
&&[H_{12},\;E_{31}]=-\frac
14\left(T+T^{-1}\right)^2E_{31}+\frac{\sf h}4\left(T-T^{-1}\right)
H_{13}^2+\frac{\sf h}2\left(T+T^{-1}\right)E_{23}H_{13}+\frac{\sf
h}{16}\left(T^3+T-T^{-1}-T^{-3}\right),\nonumber\\
&&[H_{23},\;E_{31}]=-\frac
14\left(T+T^{-1}\right)^2E_{31}+\frac{\sf h}4\left(T-T^{-1}\right)
H_{13}^2-\frac{\sf h}2\left(T+T^{-1}\right)E_{12}H_{13}+\frac{\sf
h}{16}\left(T^3+T-T^{-1}-T^{-3}\right),\nonumber\\
&&[F_{32},\;E_{21}]=F_{31}+\frac{\sf
h}4\left(T-T^{-1}\right)\left(E_{12}E_{21}+E_{23}E_{32}
\right)-\frac{\sf h}8\left(T-T^{-1}\right)H_{13}^2-\frac{\sf
h}4\left(T-T^{-1}\right)
\nonumber\\
&&\phantom{[F_3,\;E_2]=}-\frac{3{\sf
h}}{16}\left(T^2-T^{-2}\right)H_{13}-\frac{9{\sf h}}
{128}\left(T-T^{-1}\right)^3, \nonumber\\
&&[E_{12},\;E_{21}]=H_{12}+\frac
1{16}\left(T-T^{-1}\right)^2-\frac{\sf h}4\left(T-T^{-1}
\right)E_{23}E_{12},\nonumber\\
&&[E_{23},\;E_{32}]=H_{23}+\frac
1{16}\left(T-T^{-1}\right)^2+\frac{\sf h}4\left(T-T^{-1}
\right)E_{12}E_{23},\qquad [T^{\pm 1},\;E_{12}]=[T^{\pm
1},\;E_{23}]=0,\nonumber\\
&&[E_{23},\;E_{21}]=-\frac{\sf
h}4\left(T-T^{-1}\right)E_{23}^2,\qquad\qquad
[E_{12},\;E_{32}]=\frac{\sf h}4\left(T-T^{-1}\right)E_{12}^2,\\
&&[T^{\pm 1},\;E_{21}]=\mp\frac{\sf h}2\left(T^{\pm
2}+1\right)E_{23},\qquad [T^{\pm 1},\;E_{32}]=\pm\frac{\sf
h}2\left(T^{\pm 2}+1\right)E_{12},\qquad
[E_{12},\;E_{23}]=\frac1{2{\sf h}}\left(T-T^{-1}\right).\nonumber
\end{eqnarray}
The other commutators remain undeformed.
\end{prop}

\subsection{The Fermionic part}
To describe the fermionic part of the super-Jordanian quantum
superalgebra ${\cal U}_{\sf h}(sl(3|1))$, we define the elements
\begin{eqnarray}
&&H_{34}=h_{34}-\frac{{\sf h}^2}2e_{13}^2h_{13},\qquad
E_{34}=e_{34}-\frac{{\sf h}^2}4e_{13}e_{14}\left(2h_{13}+1\right),
\qquad E_{43}=e_{43},\nonumber\\
&&H_{24}=h_{24},\qquad E_{24}=e_{24},\qquad E_{42}=e_{42},\nonumber\\
&&H_{14}=h_{14}+\frac{{\sf h}^2}2e_{13}^2h_{13},\qquad
E_{14}=e_{14},\qquad E_{41}=e_{41}+\frac{{\sf
h}^2}4e_{13}e_{43}\left(2h_{13}+1\right).
\end{eqnarray}
The generators $E_{34}$, $E_{43}$, $E_{24}$, $E_{42}$, $E_{14}$
and $E_{41}$ are odd, while $H_{34}$, $H_{24}$ and $H_{14}$ are
even. The expressions (3.3), (3.4) and (3.6) constitute a
nonlinear realization of the super-Jordanian quantum superalgebra
${\cal U}_{\sf h}(sl(3|1))$ with the classical generators. Let us
just remark that the $\CC$-algebra automorphism (3.2) can be
easily extended to our construction, i.e.
\begin{equation}
\phi\left(E_{12},E_{21},H_{12},E_{23},E_{32},H_{23},E_{34},E_{43},H_{34},\cdots\right)\longrightarrow
\left(E_{23},E_{32},H_{23},E_{12},E_{12},H_{12},E_{41},E_{14},-H_{14},\cdots\right)
\end{equation}
\begin{prop}
The super-Jordanian quantum superalgebra ${\cal U}_{\sf
h}(sl(3|1))$ is then an associative superalgebra over $\CC$
spanned by
$\{H_{12},\;E_{12},\;E_{21},\;H_{23},\;E_{23},\;E_{32},\;H_{13},\;T,\;T^{-1},\;
E_{31},\;H_{34},\;E_{34},\;E_{43},\;H_{24},\;E_{24},\;E_{42},\;H_{14},$
$\;E_{14},\;E_{41}\}$, satisfying along with (3.5), the
commutation relations (we list here only the deformed commutator)
\begin{eqnarray}
&&[H_{13},\;H_{34}]=-\frac 14\left(T-T^{-1}\right)^2H_{13}, \qquad
[H_{13},\;H_{14}]=\frac
14\left(T-T^{-1}\right)^2H_{13},\nonumber\\
&& [H_{13},\;E_{14}]=\frac 12\left(T+T^{-1}\right)E_{14},\qquad
[H_{13},\;E_{43}]=\frac 12\left(T+T^{-1}\right)E_{43},\nonumber\\
&& [H_{13},\;E_{41}]=-\frac
12\left(T+T^{-1}\right)E_{41}+\frac{\sf h}2\left(T-T^{-1}\right)
E_{43}H_{13}+\frac{\sf h}4\left(T^2-T^{-2}\right)E_{43},\nonumber\\
&& [H_{13},\;E_{34}]=-\frac
12\left(T+T^{-1}\right)E_{34}-\frac{\sf h}2\left(T-T^{-1}\right)
E_{14}H_{13}-\frac{\sf h}2\left(T^2-T^{-2}\right)E_{14},\nonumber\\
&&[H_{34},\;E_{14}]=-\left(1+\frac
18\left(T-T^{-1}\right)^2\right)E_{14},\qquad
[H_{14},\;E_{43}]=\left(1+\frac
18\left(T-T^{-1}\right)^2\right)E_{43},\nonumber\\
&&[H_{34},\;E_{41}]=\left(1+\frac
18\left(T-T^{-1}\right)^2\right)E_{41}-\frac{\sf h}{16}
\left(T^2-T^{-2}\right)\left(T+T^{-1}\right)E_{43},\nonumber\\
&&[H_{34},\;E_{34}]=\frac
18\left(T-T^{-1}\right)^2E_{34}+\frac{\sf h}{16}\left(T^2-T^{-2}
\right)\left(T-T^{-1}\right)E_{14},\nonumber\\
&&[H_{34},\;E_{43}]=-\frac 18\left(T-T^{-1}\right)^2E_{43},\nonumber\\
&&[H_{34},\;T^{\pm 1}]=-\frac 14\left(T^{\pm 3}-T^{\mp
1}\right),\qquad [H_{14},\;T^{\pm 1}]=\frac 14\left(T^{\pm
3}-T^{\mp
1}\right),\nonumber\\
&&[H_{34},\;E_{31}]=\frac
14\left(T+T^{-1}\right)^2E_{31}-\frac{\sf h}4\left(T-T^{-1}\right)
H_{13}^2-\frac{\sf h}4\left(T^2-T^{-2}\right)H_{13}-\frac{\sf
h}{16}\left(T^2-T^{-2}\right)
\left(T+T^{-1}\right),\nonumber\\
&&[H_{14},\;E_{31}]=-\frac
14\left(T+T^{-1}\right)^2E_{31}+\frac{\sf h}4\left(T-T^{-1}\right)
H_{13}^2+\frac{\sf h}4\left(T^2-T^{-2}\right)H_{13}+\frac{\sf
h}{16}\left(T^2-T^{-2}\right)
\left(T+T^{-1}\right),\nonumber\\
&&[H_{14},\;E_{34}]=-\left(1+\frac
18\left(T-T^{-1}\right)^2\right)E_{34}-\frac{\sf h}{16}
\left(T^2-T^{-2}\right)\left(T+T^{-1}\right)E_{43}, \nonumber\\
&&[T^{\pm 1},\;E_{34}]=\pm \frac{\sf h}2\left(T^{\pm
2}+1\right)E_{14},\qquad [T^{\pm 1},\;E_{41}]=\mp\frac{\sf
h}2\left(T^{\pm
2}+1\right)E_{43},\nonumber\\
&& [E_{43},\;E_{14}]=\frac 1{2{\sf h}}\left(T-T^{-1}\right),\nonumber\\
&&[E_{34},\;E_{43}]=H_{34}-\frac
1{16}\left(T-T^{-1}\right)^2-\frac{\sf h}4\left(T-T^{-1}
\right)E_{14}E_{43},\nonumber\\
&&[E_{14},\;E_{41}]=H_{14}+\frac
1{16}\left(T-T^{-1}\right)^2+\frac{\sf h}4\left(T-T^{-1}
\right)E_{43}E_{14},\nonumber\\
&&[E_{43},\;E_{31}]=\frac{\sf
h}4\left(T-T^{-1}\right)E_{34}+\frac{\sf h}2\left(T-T^{-1}\right)
E_{14}E_{31}-\frac{{\sf h}^2}4E_{14}H_{13}^2-\frac{3{\sf
h}^2}8\left(T+T^{-1}\right)E_{14}H_{13}
\nonumber\\
&&\phantom{[E_2,\;F_1]=}-\frac{{\sf h}^2}2E_{14}-\frac{15{\sf
h}^2}{64}\left(T-T^{-1}\right)^2
E_{14},\nonumber\\
&&[E_{41},\;E_{31}]=\frac{\sf
h}4\left(T-T^{-1}\right)E_{41}-\frac{\sf h}2\left(T-T^{-1}\right)
E_{43}E_{31}+\frac{{\sf h}^2}4E_{43}H_{13}^2+\frac{3{\sf
h}^2}8\left(T+T^{-1}\right)E_{43}H_{13}
\nonumber\\
&&\phantom{[E_2,\;F_1]=}+\frac{{\sf h}^2}2E_{43}+\frac{15{\sf
h}^2}{64}\left(T-T^{-1}\right)^2
E_{43},\nonumber\\
&&[E_{43},\;E_{32}]=F_{42}+\frac{\sf
h}4\left(T-T^{-1}\right)E_{12}E_{43},\nonumber\\
&&E_{34}^2=\frac{\sf h}4\left(T-T^{-1}\right) E_{14}E_{34},\qquad
E_{41}^2=-\frac{\sf
h}4\left(T-T^{-1}\right)E_{43}E_{41},\nonumber\\
&&[T^{\pm 1},\;E_{14}]=0,\qquad [T^{\pm 1},\;E_{43}]=0, \qquad
[T^{\pm 1},\;E_{24}]=0,\qquad
[T^{\pm 1},\;E_{42}]=0,\nonumber\\
&&[E_{34},\;E_{41}]=F_{31}-\frac{\sf
h}4\left(T-T^{-1}\right)E_{43}E_{34}+\frac{\sf h}4
\left(T-T^{-1}\right)E_{14}F_{41}-\frac{\sf
h}8\left(T-T^{-1}\right)H_{13}^2-\frac{\sf h}8
\left(T^2-T^{-2}\right)H_{13}^2\nonumber\\
&&\phantom{[F_3,\;E_2]=}-\frac{\sf
h}{16}H_{13}\left(T^2-T^{-2}\right)+ \frac{7{\sf
h}}{128}\left(T-T^{-1}\right)^3.
\end{eqnarray}
\end{prop}

The coalgebraic structure will be presented, for the general case,
in the following section.

\section{${\cal U}(sl(N|1))$: Generalization}
\setcounter{equation}{0}

From the above studies, it is easy to see that:
\begin{prop}
The superalgebra ${\cal U}_{\sf h}(sl(N|1))$ can be realized via
the nonlinear map:
\begin{eqnarray}
&&T^{\pm 1}=\pm {\sf h}e_{1N}+\sqrt{1+{\sf h}^2e_{1N}^2},\qquad
H_{1N}=\sqrt{1+{\sf h}^2e_{1N}^2}h_{1N},\qquad
E_{N1}=e_{N1}-\frac{{\sf
h}^2}4e_{1N}\left(h_{1N}^2-1\right),\nonumber\\
&&H_{ij}=h_{ij}+\frac{{\sf
h}^2}2\left(\delta_{i1}+\delta_{jN}\right)e_{1N}^2h_{1N},
\qquad\qquad
i<j\in\{1,2,\cdots,N\}\;\;\hbox{and}\;\;(i,j)\neq(1,N),\nonumber\\
&&E_{ij}=e_{ij},\qquad\qquad
i<j\in\{1,2,\cdots,N\}\;\;\hbox{and}\;\;(i,j)\neq(1,N),\nonumber\\
&&E_{ji}=e_{ji}+\frac{{\sf
h}^2}4\left(\delta_{i1}e_{jN}-\delta_{Nj}e_{1i}\right)
\left(2h_{1N}+1\right),\qquad\qquad
i<j\in\{1,2,\cdots,N\}\;\;\hbox{and}\;\;(i,j)\neq(1,N),\nonumber\\
&&H_{i,N+1}=h_{i,N+1}+\frac{{\sf
h}^2}2\left(\delta_{i1}-\delta_{iN}\right)e_{1N}^2h_{1N},
\qquad\qquad i\in\{1,2,\cdots,N\},\nonumber\\
&&E_{i,N+1}=e_{i,N+1}-\frac{{\sf
h}^2}4\delta_{iN}e_{1,N+1}e_{1N}\left(2h_{1N}+1\right),
\qquad\qquad i\in\{1,2,\cdots,N\},\nonumber\\
&&E_{N+1,i}=e_{N+1,i}+\frac{{\sf
h}^2}4\delta_{iN}e_{N+1,N}e_{1N}\left(2h_{1N}+1\right),
\qquad\qquad i\in\{1,2,\cdots,N\},
\end{eqnarray}
with the coproducts
\begin{eqnarray}
&&\Delta(H_{1N})=H_{1N}\otimes T+T^{-1}\otimes H_{1N},\qquad
\Delta(T^{\pm 1})=T^{\pm 1}\otimes T^{\pm 1},\qquad
\Delta(E_{N1})=E_{N1}\otimes T+T^{-1}\otimes E_{N1},\nonumber\\
&&\Delta\left(H_{ij}\right)=H_{ij}\otimes 1+1\otimes H_{ij}-\frac
14\left(\delta_{i1}+
\delta_{jN}\right)\biggl(TH_{1N}\otimes\left(1-T^2\right)+\left(1-T^{-2}\right)\otimes
T^{-1}H_{1N}\biggr),\nonumber\\
&&\phantom{xxxxxxxxxx}i<j\in\{1,2,\cdots,N\}\;\;\hbox{and}\;\;(i,j)\neq(1,N),\nonumber\\
&&\Delta\left(E_{ij}\right)=E_{ij}\otimes
T^{-\left(\delta_{i1}+\delta_{jN}\right)/2}+
T^{\left(\delta_{i1}+\delta_{jN}\right)/2}\otimes E_{ij},\qquad
i<j\in\{1,2,\cdots,N\}\;\;\hbox{and}\;\;(i,j)\neq(1,N),\nonumber\\
&&\Delta\left(E_{ji}\right)=E_{ji}\otimes
T^{\left(\delta_{i1}+\delta_{jN}\right)/2}+
T^{-\left(\delta_{i1}+\delta_{jN}\right)/2}\otimes
E_{ji}+\frac{\sf h}4T^{-1}\left(-\delta_{i1}
E_{jN}+\delta_{jN}E_{1i}\right)\otimes
\left(T^{-1/2}H_{1N}+\right.\nonumber\\
&&\phantom{\Delta\left(E_{ji}\right)=}\left.
H_{1N}T^{-1/2}\right)-\frac{\sf h}4\left(T^{1/2}
H_{1N}+H_{1N}T^{1/2}\right)\otimes
T\left(-\delta_{i1}E_{jN}+\delta_{jN}E_{1i}\right)\nonumber\\
&& \phantom{xxxxxxxxxx}
i<j\in\{1,2,\cdots,N\}\;\;\hbox{and}\;\;(i,j)\neq(1,N),\nonumber\\
&&\Delta\left(H_{i,N+1}\right)=H_{i,N+1}\otimes 1+1\otimes
H_{i,N+1}-\frac 14\left(\delta_{i1}-
\delta_{iN}\right)\biggl(TH_{1N}\otimes\left(1-T^2\right)+\left(1-T^{-2}\right)\otimes
T^{-1}H_{1N}\biggr),\nonumber\\
&&\phantom{xxxxxxxxxx}i\in\{1,2,\cdots,N\},\nonumber\\
&&\Delta\left(E_{i,N+1}\right)=E_{i,N+1}\otimes
T^{\left(-\delta_{i1}+\delta_{iN}\right)/2}+
T^{\left(\delta_{i1}-\delta_{iN}\right)/2}\otimes E_{ji}+\frac{\sf
h}4T^{-1}\delta_{iN}
E_{1N}\otimes \left(T^{-1/2}H_{1N}+\right.\nonumber\\
&&\phantom{\Delta\left(E_{i,N+1}\right)=}\left.
H_{1N}T^{-1/2}\right)-\frac{\sf h}4
\left(T^{1/2}H_{1N}+H_{1N}T^{1/2}\right)\otimes
\delta_{iN}TE_{1N}\qquad
i\in\{1,2,\cdots,N\},\nonumber\\
&&\Delta\left(E_{N+1,i}\right)=E_{N+1,i}\otimes
T^{\left(\delta_{i1}-\delta_{iN}\right)/2}+
T^{\left(-\delta_{i1}+\delta_{iN}\right)/2}\otimes
E_{N+1,i}-\frac{\sf h}4T^{-1}\delta_{i1}
E_{N+1,N}\otimes \left(T^{-1/2}H_{1N}+\right.\nonumber\\
&&\phantom{\Delta\left(E_{i,N+1}\right)=}\left.
H_{1N}T^{-1/2}\right)+\frac{\sf h}4
\left(T^{1/2}H_{1N}+H_{1N}T^{1/2}\right)\otimes
\delta_{i1}TE_{N+1,N}\qquad i\in\{1,2,\cdots,N\}.
\end{eqnarray}
\end{prop}

The commutator rules of ${\cal U}_{\sf h}(sl(N|1))$ can be
evaluated by direct calculations.

Paralleling the earlier cases, the universal ${\cal R}_{\sf
h}$-matrix of ${\cal U}_{\sf h}(sl(N|1))$ is given by
\begin{eqnarray}
&& {\cal R}_{\sf h}=\exp\biggl(-{\sf h}E_{1N}\otimes
TH_{1N}\biggr)\exp\biggl({\sf h}TH_{1N} \otimes E_{1N}\biggr),
\end{eqnarray}
where $E_{1N}={\sf h}^{-1}\ln T$. This element can be connected to
the results obtained by the contraction process by a suitable
twist operator that can be derived as a series expansion in ${\sf
h}$.

\section{Conclusion}
\setcounter{equation}{0}

In general, a class of nonlinear maps exists relating the
Jordanian quantum (super)algebras and their classical analogues.
Here we have used a particular map realizing Jordanian ${\cal
U}_{\sf h}(sl(N|1))$ superalgebra for an arbitrary $N$. This map
arises naturally form our contraction process defined in (1.1).
Let us just recall that the more important advantages of our
procedure are:

\begin{itemize}
\item The algebraic commutation relations are deformed.
\item The coalgebraic structure is simpler.
\item The map obtained, by  our contraction process, permits immediate
explicit construction of the finite-dimensional irreps.
\item The Ohn's ${\cal U}_{\sf h}(sl(2))$ algebra is embedded as a Hopf
subalgebra in our construction. Therefore, the Jordanian ${\cal
U}_{\sf h}(sl(N|1))$ superalgebra arising from our method
corresponds to the classical $r$-matrix $r=h_{1N}\otimes
e_{1N}-e_{1N}\otimes h_{1N}$.
\end{itemize}

\vskip 0.5cm

\noindent{\bf Acknowledgments:} One of us (BA) wants to thank
Jean-Loup Roussel for precious help. This work is supported by a
grant of "La Fondation Charles de Gaulle". The work of R.C. is
partially supported by a grant of DAE(BRNS), Government of India.

\end{document}